\documentclass[journal,12pt,onecolumn,draftclsnofoot]{IEEEtran}

\usepackage{ifpdf}

\usepackage[ruled]{algorithm2e}
\usepackage{color}
\usepackage{url}
\ifpdf
\usepackage[pdftex]{graphicx}
\else
\usepackage[dvips]{graphicx}
\fi
\usepackage[caption=false,font=footnotesize]{subfig}
\usepackage{comment}
\usepackage{amsmath,amsfonts} 
\usepackage[noadjust]{cite}
\usepackage{nomencl}
\newtheorem{theorem}{\textbf{Theorem}}

\newtheorem{lemma}{\textbf{Lemma}}

\usepackage{xcolor}

\definecolor{ggray}{gray}{.9}

\newcommand{\ud}[1]{_\mathrm{#1}}
\newcommand{\up}[1]{^\mathrm{#1}}
\newcommand\Tstrut{\rule{0pt}{2.6ex}}         
\newcommand\Bstrut{\rule[-0.9ex]{0pt}{0pt}}   

\newcommand{\highlight}[1]{\colorbox{pink}{#1}}
\newcommand{\compare}[1]{\colorbox{ggray}{#1}}

\IEEEoverridecommandlockouts
\begin{document}

\title{Optimal Regulation Response of Batteries \\Under  Cycle Aging Mechanisms}

\author{
    Bolun Xu,
    Yuanyuan Shi,
    Daniel S. Kirschen
    and Baosen Zhang
\thanks{The authors are with the Department of Electrical Engineering,
University of Washington,
Seattle, Washington 98125,
(e-mail:\{xubolun, yyshi, kirschen, zhangbao\}@uw.edu). This work has been supported in part by the University of Washington Clean Energy Institute. }
}

\maketitle
\makenomenclature

\begin{abstract}
When providing frequency regulation in a pay-for-performance market, batteries need to carefully balance the trade-off between following regulation signals and their degradation costs in real-time. Existing battery control strategies either do not consider mismatch penalties in pay-for-performance markets, or cannot accurately account for battery cycle aging mechanism during operation. This paper derives an online control policy that minimizes a battery owner's operating cost for providing frequency regulation in a pay-for-performance market. The proposed policy considers an accurate electrochemical battery cycle aging model, and is applicable to most types of battery cells. It has a threshold structure, and achieves near-optimal performance with respect to an offline controller that has complete future information. We explicitly characterize this gap and show it is independent of the duration of operation. Simulation results with both synthetic and real regulation traces are conducted to illustrate the theoretical results.
\end{abstract}


\IEEEpeerreviewmaketitle

\section{Introduction}

Batteries are becoming a key provider of frequency regulation, a power system ancillary service needed to maintain the system frequency.
Following Federal Energy Regulatory Commission (FERC) Order 755~\cite{xu2016comparison} in 2011, most major U.S. system operators have implemented pay-for-performance regulation markets. In these markets, a participant's payment for providing the regulation service depends not only on the regulation capacity it provides, but also on the  accuracy of its  response to the regulation instruction. The instruction takes the form of a signal that is sent every 2 to 6 seconds, representing the amount of active power a participant should inject or withdraw. A participant is penalized if it deviates from the received regulation signal. Since batteries have much faster ramp rates compared to traditional generators, they are among the most competitive providers in these fast regulation markets~(in 2016 41\% of regulation in PJM was provided by batteries~\cite{pjm_update}). The importance of batteries is likely to increase further as the size of the fast regulation markets grows in response to the growing penetration of renewable generation.

Although batteries can achieve near perfect accuracy in the provision of regulation ~\cite{pjm_update}, it is not always clear that a battery should exactly follow the regulation signal to maximize its gain from participating in the regulation markets. The optimal action of a battery should balance the profit from providing regulation with its operating cost, which is mainly driven by the degradation caused by the charge and discharge cycles~\cite{dunn2011electrical,zakeri2015electrical}. In particular, deep cycles--charging/discharging cycles that use a significant amount of active materials--tend to dramatically reduce battery life~\cite{vetter2005ageing}. Indeed, a naive battery controller that attempts to follow regulation signals without considering cycle degradation could destroy a battery in a matter of months. Instead, a controller should strategically choose to suffer some performance penalty to avoid deep cycles.  However, designing a better controller is not a trivial task, since it is difficult to tell whether the battery is undergoing a deep or a shallow cycle without future knowledge. Most commercial controllers sidestep this issue by setting hard limits on the battery state of charge, which can limit the profitability of batteries and artificially increase the need for more regulation resources.

In this paper, we overcome the challenge of balancing regulation performance and reducing battery degradation by designing an online control policy that is nearly optimal:
\begin{enumerate}
	\item It achieves a bounded optimality gap compared with an optimal offline policy that has full information about future regulation signals.
	\item This gap is independent of the duration of operationand can be characterized exactly.
\end{enumerate}
The key to this control policy is a more thorough \emph{algorithmic understanding} of the battery aging process with respect to charge/discharge cycles.

The rest of this paper is organized as follows.
Section~\ref{sec:lit} introduces the battery degradation mechanism and compares the proposed policy with similar works. Section~\ref{sec:model} introduces the model and problem formulations. Section~\ref{sec:policy} describes the proposed control strategy, while the optimality proof is given in the appendix. The performance of this strategy is validated through simulations in Section~\ref{sec:simulation}, and Section~\ref{sec:con} concludes the paper.


\section{Background and Literature Review}\label{sec:lit}

Two key operational factors that accelerates battery degradation are \emph{high current rates} and \emph{deep cycles}. Since system operators normally require battery storage to have at least a 15 minutes capacity to provide regulation~\cite{xu2016comparison}, the effect of current rate on degradation is relatively small~\cite{wang2014degradation}, and the cycle aging effect dominates in regulation markets.

Cycle aging depends nonlinearly on the \emph{charge~/~discharging cycle depths} in most static electrochemical batteries~\cite{ruetschi2004aging, byrne2012estimating, xu2016modeling, wang2014degradation,millner2010modeling,laresgoiti2015modeling}. For example, a 7~Wh Lithium Nickel Manganese Cobalt Oxide (NMC) battery provides 35~kWh lifetime energy throughput if cycled at 10\% depth, while only 3.5~kWh at 100\% cycle depth~\cite{ecker2014calendar}. Physically, a cycle is a fatigue event in which the battery first deviates and later returns to a certain state of charge (SoC) level.  The depth is independent of the current rate, or the starting or ending point of this cycle. It is thus similar to a material fatigue process.

\begin{figure}[t]%
	\centering
	\subfloat[Regulation instruction for consumption and injection.]{
		\includegraphics[trim = 0mm 0mm 0mm 0mm, clip, width = .7\columnwidth]{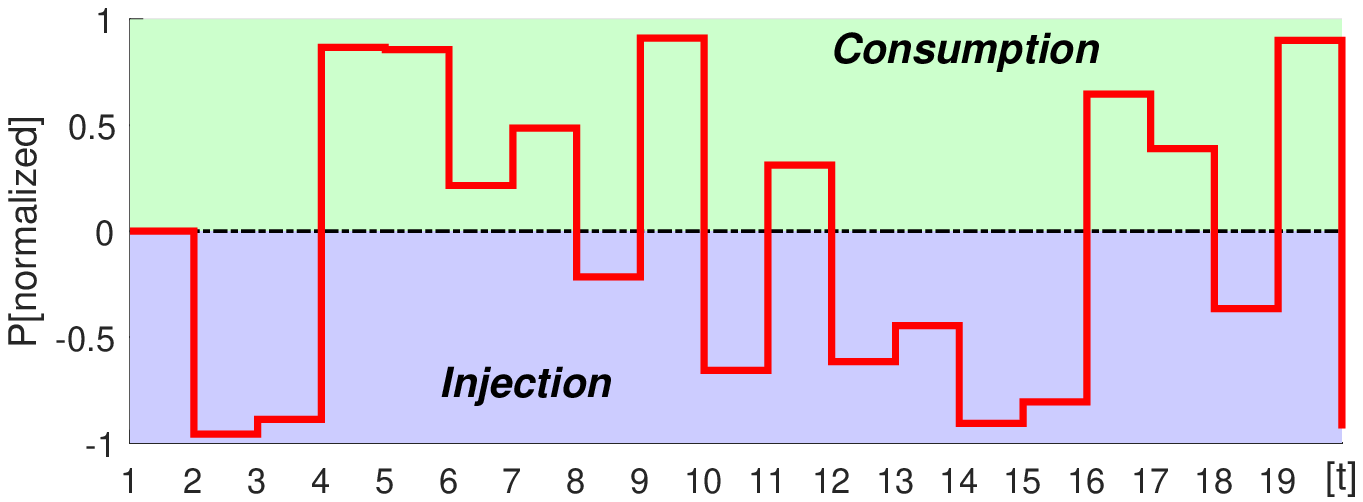}
		\label{fig:depth1}%
	}
	\\
	\subfloat[SoC of a battery that follows the regulation instruction.]{
		\includegraphics[trim = 0mm 0mm 0mm 0mm, clip, width = .7\columnwidth]{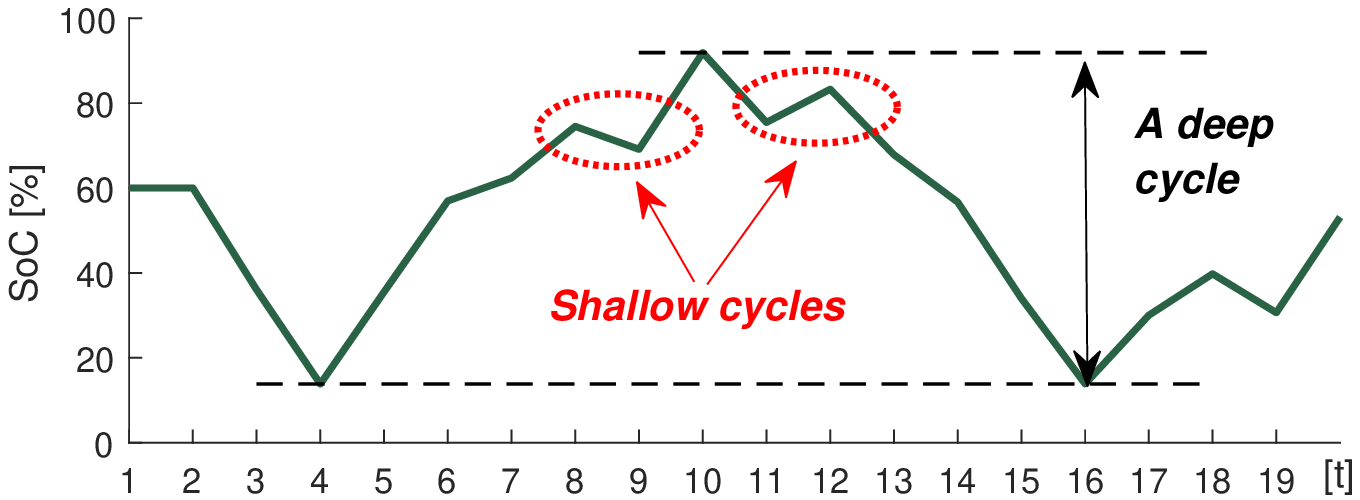}
		\label{fig:depth2}%
	}
	\caption{\footnotesize An example of a battery following regulation instructions. A deep cycle and two shallow cycles are contained in the operation between $t=4$ and $t=16$.}%
	\label{fig:depth}
\end{figure}

Fig~\ref{fig:depth} shows an example battery operation and the resulting SoC profile. A deep cycle can be clearly observed between $t=4$ and $t=16$. This deep cycle is interrupted by two shallower cycles. In practice, a systematic cycle counting process must be adopted to account for arbitrary SoC profiles. Currently, the most popular method is the \emph{rainflow} cycle counting algorithm. This algorithm is used extensively in materials fatigue stress analysis to count cycles and quantify their depth, and has also been extensively applied to battery life assessment~\cite{xu2016modeling, muenzel2015multi}.

The main technical challenge with the rainflow algorithm is that it does not have an analytical expression that can be used directly in optimization problems~\cite{amzallag1994standardization} and much of the focus has been on assessing the degradation of batteries \emph{after the fact}.
Several efforts have been made to incorporate simplified rainflow methods into the optimization of battery operation through model predictive control~\cite{koller2013defining,fortenbacher2017modeling}, dynamic programming~\cite{amzallag1994standardization}, or stochastic programming~\cite{he2015optimal}. These approaches all make additional simplifications to the cycle aging model, but are still too complex for controlling the response to a regulation in real-time every 2 to 6 seconds.

{The control policy proposed in this paper uses a much simpler approach which keeps the battery's SoC between an upper and a lower bound to avoid deep cycles. This policy reacts to new regulation instructions instantaneously and achieves near-optimal control without having to solve an optimization problem in real-time. The idea of regulating SoC has been extensively incorporated in previous battery control strategies.} Common approaches include constraining the SoC within pre-fixed bounds~\cite{oudalov2007optimizing}, or maintaining the SoC at a target level with a proportional-integral (PI) controller~\cite{borsche2013power}. These approaches are simple to implement and are effective at alleviating battery aging~\cite{xu2014bess}, but their heuristic settings are not responsive to market prices and thus limit the profitability of the batteries.

In contrast, we derive the SoC thresholds from first principles based on the regulation market prices and the battery cycle aging mechanism. A related work ~\cite{zhang2016profit} proposed a profit-maximizing battery control strategy for frequency control based on price signals, but the cost function is too simple to model the cycle-based battery aging mechanism. Our policy incorporates the battery aging model into the SoC threshold calculations, improving model accuracy and making it applicable to most electrochemical battery cells. This battery control policy also applies to any market application that has stochastic dispatch but constant prices over a specific period, such as behind-the-meter peak shaving~\cite{shi2016leveraging}, or improving the penetration of renewable generations~\cite{bitar2011role}.

\begin{figure*}[t]%
	\centering
	\subfloat[The full example profile]{
		\includegraphics[trim = 0mm 0mm 0mm 0mm, clip, width = .3\columnwidth]{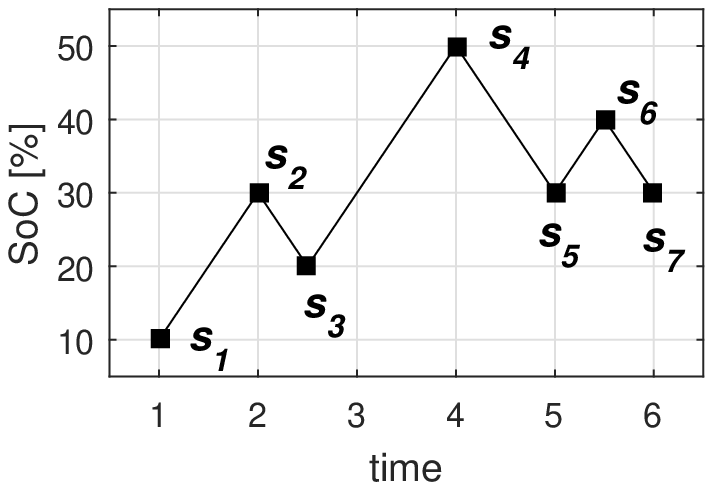}
		\label{fig:rf5}%
	}
	\subfloat[Extracted full cycles]{
		\includegraphics[trim = 0mm 0mm 0mm 0mm, clip, width = .3\columnwidth]{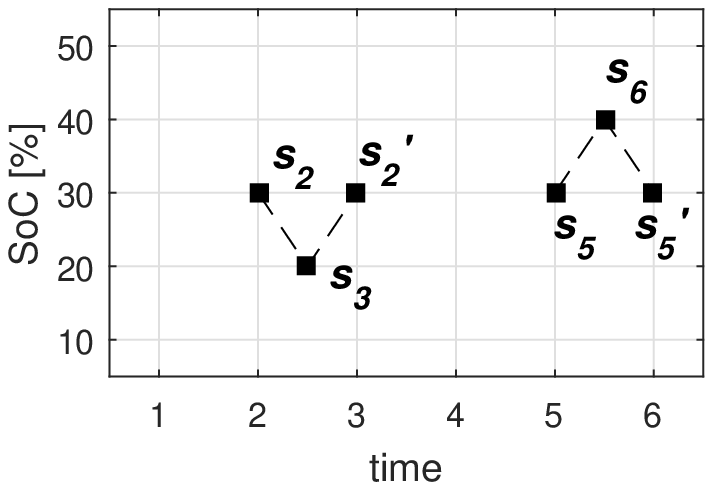}
		\label{fig:rf6}%
	}
	\subfloat[The remaining residue profile]{
		\includegraphics[trim = 0mm 0mm 0mm 0mm, clip, width = .3\columnwidth]{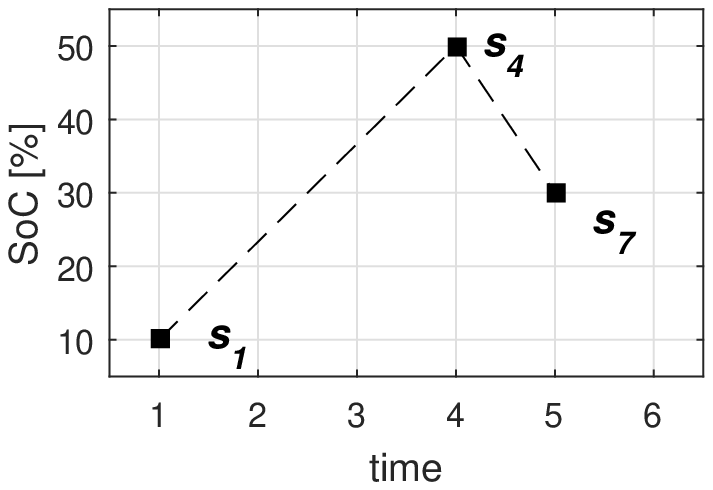}
		\label{fig:rf7}%
	}
	\caption{\footnotesize Rainflow cycle counting result example.}%
	\label{Fig:rf}
\end{figure*}

\section{Models for Battery and Regulation}~\label{sec:model}
In this section we describe the battery model and how the rainflow algorithm can be used to calculate the cost of battery cycle aging from control actions. Then we model the regulation market performance settlement process, and state the main optimization problem of the paper: \emph{how to balance profit from regulation and the degradation cost of battery operation in an online fashion.}

\subsection{Battery Operations}
We consider an operation defined over finite discrete control time intervals $n\in \{1,\dotsc, N\}$, $N\in \mathbb{N}$. Each control time interval has a duration of $T$, and the entire operation has a duration of $TN$. Let $e_n$ be the energy stored in the battery--the state of charge (SoC)--at time step $n$. The battery can either charge an amount $c_n$  or discharge an amount $d_n$ during each interval. Its state of charge is given by the following linear difference equation:
\begin{align}
    e_{n+1} - e_{n} =& T\eta\ud{c} c_{n} - (T/\eta\ud{d}) d_n\label{eq:bat_soc}.
\end{align}
For a given battery, it has the following operation constraints
\begin{align}
    \underline{e} \leq e_n &\leq \overline{e}
    \label{eq:bat_con1}\\
    0 \leq c_{n} &\leq P
    \label{eq:bat_con2}\\
    0 \leq d_n &\leq P
    \label{eq:bat_con3}\\
    c_n = 0 \mbox{ or } & d_n=0 \label{eq:bat_con4}
\end{align}
where $\underline{e}$ and $\overline{e}$ are the minimum and maximum energy capacity of the battery, $\eta\ud{c}$ and $\eta\ud{d}$ are the charge and discharge efficiency, $P$ is the battery power rating, and \eqref{eq:bat_con4} avoids trivial solutions by preventing the battery from charging and discharging at the same time.


\subsection{Cycle Counting via Rainflow}
The cycle aging of electrochemical battery cells is evaluated based on stress cycles, and the rainflow method identify cycles from local extrema in a SoC profile. A local extrema point indicates that the battery switched from charge mode to discharge mode, or vice versa. Consider a SoC profile with local extrema $s_1, s_2, \dotsc$ (as in Fig.~\ref{fig:rf5}), the rainflow method identifies cycles according to the following procedure~\cite{ramadesigan2012modeling}
\begin{enumerate}
    \item Start from the beginning of the profile.
    \item Calculate $\Delta s_1 = |s_1-s_2|$, $\Delta s_2 = |s_2-s_3|$, $\Delta s_3 = |s_3-s_4|$.
    \item If $\Delta s_2\leq \Delta s_1$ and $\Delta s_2 \leq \Delta s_3$, then a full cycle of depth $\Delta s_2$ associated with $s_2$ and $s_3$ has been identified. Remove $s_2$ and $s_3 $ from the profile, and repeat the identification using points $s_1$, $s_2$, $s_5$, $s_6$...
    \item If a cycle has not been identified, shift the identification forward and repeat the identification using points $s_2$, $s_3$, $s_4$, $s_5$...
    \item The identification is repeated until no more full cycles can be identified throughout the remaining profile.
\end{enumerate}
Fig.~\ref{fig:rf6} shows full cycles that are extracted from the stress profile in Fig.~\ref{fig:rf5}, and Fig.~\ref{fig:rf7} is the remaining rainflow residue profile that contains no full cycles. To classify the aging stress caused by the residue profile, we employ the half cycle method for rainflow residue treatment~\cite{marsh2016review}. In this method, the rainflow residue is decomposed into half cycles, each half cycle causes half the aging stress of a full cycle of the same depth. A half cycle is between each two adjacent local extrema in the rainflow residue. We define a half cycle with increasing SoC as a charge half cycle, and a half cycle with decreasing SoC as a discharge half cycle. The residue profile in Fig.~\ref{fig:rf7} contains a charge half cycle of depth 40\%, and a discharge half cycle of depth 20\%.

Let $\mathbf{c}$ and $\mathbf{d}$ represents a series of charge and discharge controls over a time period, then cycles can be calculated directly from $(\mathbf{c}, \mathbf{d})$ based on the battery energy capacity $E$~(regardless of the initial starting point $e_0$). Let $\mathrm{Rainflow}$ be the rainflow counting algorithm, then
\begin{align}\label{eq:rf}
    (\mathbf{u}, \mathbf{v}, \mathbf{w}) = \mathrm{Rainflow}\Big(\frac{T\eta\ud{c}}{E} \mathbf{c}- \frac{T}{\eta\ud{d}E} \mathbf{d}\Big)
\end{align}
where $\mathbf{u}$ is the set of all full cycle depths, $\mathbf{v}$ for charge half cycles, and $\mathbf{w}$ for discharge half cycles. 

\subsection{Battery Degradation Cost}
The cycle aging model employs a cycle depth stress function $\Phi(u):[0,1]\to \mathbb{R}^+$ to model the life loss from a single cycle of depth $u$. This function indicates that if a battery cell is repetitively cycled with depth $u$, then it can operate $1/\Phi(u)$ number of cycles before reaching its end of life. $\Phi(u)$ is a convex function for most types of electrochemical batteries~\cite{millner2010modeling,ecker2014calendar,ruetschi2004aging, byrne2012estimating, wang2014degradation}, an example polynomial form of $\Phi(u)$ is as $\alpha u^{\beta}$~\cite{laresgoiti2015modeling}. Because cycle aging is a cumulative fatigue process~\cite{millner2010modeling,ecker2014calendar}, the total life loss $\Delta L$ is the sum of the life loss from all cycles
\begin{align}\label{eq:deltaL}
   \Delta L(\mathbf{u}, \mathbf{v}, \mathbf{w}) = \sum_{i=1}^{|\mathbf{u}|} \Phi(u_i)+\sum_{i=1}^{|\mathbf{v}|} \frac{\Phi(v_i)}{2}+\sum_{i=1}^{|\mathbf{w}|} \frac{\Phi(w_i)}{2}
\end{align}
where $|\mathbf{u}|$ is the cardinality of $\mathbf{u}$. For example, to calculate cycle aging for the profile in Fig.~\ref{fig:rf5}, we set $u_1 = 0.1$, $u_2 = 0.1$, $v_1 = 0.4$ and $w_1 = 0.2$. If we subsitute the rainflow algorithm as in \eqref{eq:rf} into \eqref{eq:deltaL}, the incremental cycle aging can therefore be written as a function of the control actions $\mathbf{c}$ and $\mathbf{d}$. Let $R$ be the battery cell replacement price in \$/MWh and hence $ER$ is the replacement cost. The cycle aging cost function $J\ud{cyc}(\mathbf{c}, \mathbf{d})$ is
\begin{align} \label{eqn:Jcyc}
    J\ud{cyc}(\mathbf{c}, \mathbf{d}) = \Delta L(\mathbf{c}, \mathbf{d})ER\,.
\end{align}

\subsection{Market Settlement Model}

The system operator clears the regulation market to determine the performance penalty price before a dispatch interval~\cite{xu2016comparison}.
Instead of going into cumbersome market details, we generalize the ex-post regulation market settlement model from the perspective of a participant as follows. We assume market prices are constant throughout the operation period.  Suppose a regulation participant pays a constant positive over-response price $\theta\in R^+$ ($\$/$MWh) for surplus injections or deficient demands during each dispatch interval, and a constant under-response price $\pi\in R^+$ ($\$/$MWh) for deficient injections or surplus demands. The performance penalty model $J\ud{reg}$ calculates the market settlement cost for performance the regulation
\begin{align}
    J\ud{reg}(\mathbf{c}, \mathbf{d}) = &T\theta\textstyle\sum_{n=1}^N|c_n - d_n - r_n|^+ \nonumber\\
    &+T\pi\textstyle\sum_{n=1}^N|r_n-c_n+d_n|^+\,,
\end{align}
where $r_n\in [-P, P]$ is the instructed regulation dispatch set-point for the dispatch interval $n$, with the convention positive values in $r_n$ represents charge instructions.

\subsection{Optimization Problem}
If the regulation instruction $\mathbf{r}$ is known, then the optimization problem is:
\begin{subequations} \label{eqn:main}
  \begin{align}
    \min_{\mathbf{c}, \mathbf{d}} \; &J(\mathbf{c}, \mathbf{d}, \mathbf{r}):= J\ud{cyc}(\mathbf{c}, \mathbf{d}) + J\ud{reg}(\mathbf{c}, \mathbf{d}, \mathbf{r}) \\
    \mbox{s.t. } &\eqref{eq:bat_con1}-\eqref{eq:bat_con4}.
  \end{align}
\end{subequations}
However, this problem is inherently online: the charging and discharging decisions must be made at each time step without knowing the future realization of the regulation instruction $\mathbf{r}$. Therefore we seek an \emph{online policy} that will determine $c_n$ and $d_n$ at time step $n$ with only past information. Note we do not assume any information about the future realization of the regulation signal is known~(e.g., it need not come from a stochastic process).

\section{Proposed Online Control Policy}\label{sec:policy}

\begin{figure}[t]%
	\centering
	\subfloat[Regulation instruction vs. response]{
		\includegraphics[trim = 0mm 0mm 0mm 0mm, clip, width = .7\columnwidth]{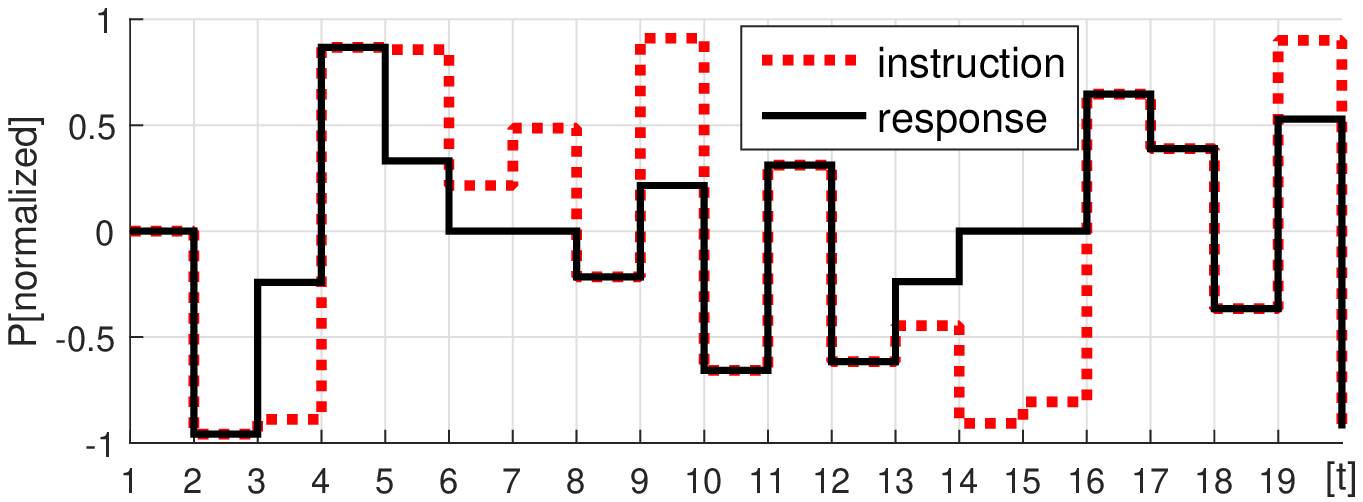}
		\label{fig:con1}%
	}
	\\
	\subfloat[Controlled vs. uncontrolled SoC profile]{
		\includegraphics[trim = 0mm 0mm 0mm 0mm, clip, width = .7\columnwidth]{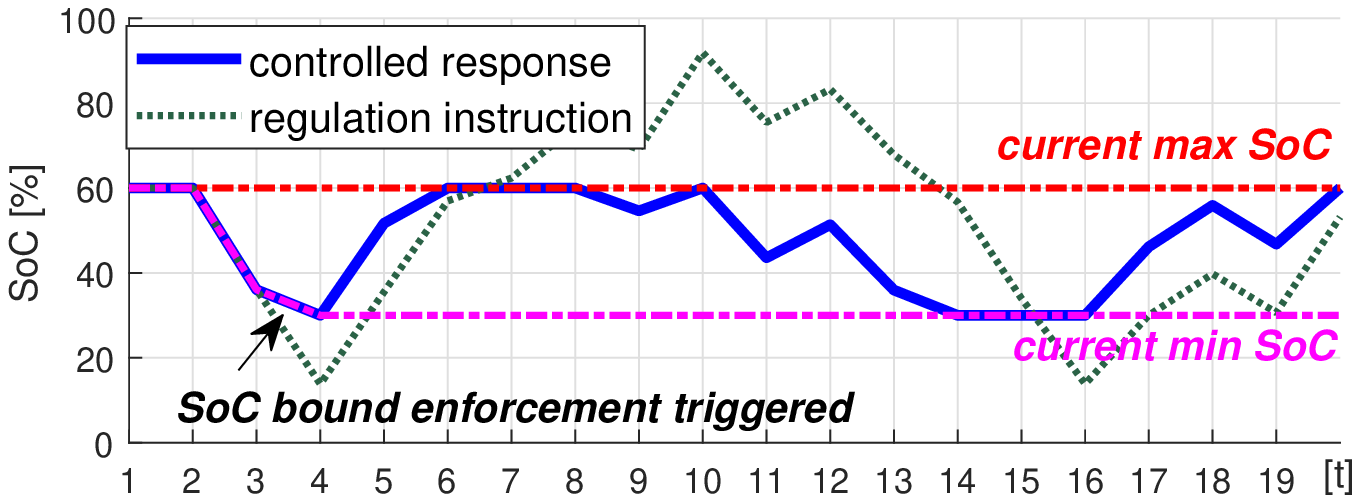}
		\label{fig:con2}%
	}
	\caption{\footnotesize Proposed control policy illustration. The policy keeps track the current maximum and minimum SoC level. When the distance in between reaches the calculated threshold $\hat{u}$, the policy starts to constrain the response. 
	}%
	\label{fig:con}
\end{figure}
We propose an online battery control policy that balances the cost of deviating from the regulation signal and the cycle aging cost of batteries while satisfying operation constraints. This policy takes a \emph{threshold form} and achieves an optimality gap that is \emph{independent of the total number of time steps}. Therefore in term of regret, this policy achieves the strongest possible result: the regret do not grow with time. Note we assume the regulation capacity has already been fixed in the previous capacity settlement stage~\cite{xu2014bess}.

\subsection{Control Policy Formulation}
The key part of the control policy is to calculate thresholds that bounds the SoC of the battery as functions of the deviation penalty and degradation cost. Let $\hat{u}$ denote this bound on the SoC, and it is given by:
\begin{align}\label{eq:pol3}
    \hat{u} = \dot{\Phi}^{-1}\Big(\frac{\pi\eta\ud{d}+ \theta/\eta\ud{c}}{R}\Big) 
\end{align}
where $\dot{\Phi}^{-1}(\cdot)$ is the inverse function of the derivative of the cycle stress function $\Phi(\cdot)$.

\begin{algorithm}[!htb]\label{alg}
\SetAlgoLined
\KwResult{Determine battery dispatch point $c_n$, $d_n$}
 \tcp{initialization}
 set $\Phi\Big(\frac{\pi\eta\ud{d}+ \theta/\eta\ud{c}}{R}\Big) \to \hat{u}$, $e_0 \to e\up{max}_{0}$, $e_0 \to e\up{min}_{0}$\;
 \While{$n\leq N$}{
  \tcp{read $e_n$ and update controller}
  set $\max\{e\up{max}_{n-1}, e_n\} \to e\up{max}_n$, $\min\{e\up{min}_{n-1}, e_n\} \to e\up{max}_n$\;
  set $\min\{\overline{e}, e\up{min}_n + \hat{u}E\} \to \overline{e}^g_n$\;
  set $\max\{\underline{e}, e\up{max}_n + \hat{u}E\} \to \underline{e}^g_n$\;
  \tcp{read $r_n$ and enforce soc bound}
  \eIf{$r_n \geq 0$}{
   set $\min\Big\{\frac{1}{T\eta\ud{c}}(\overline{e}^g-e_n), r_n\Big\} \to c_n$, $0 \to d_n$ \;
   }{
   set $0\to c_n$, $\min\Big\{\frac{\eta\ud{d}}{T}(e_n-\underline{e}^g), r_n\Big\}\to d_n$ \;
  }
  \tcp{wait until next control interval}
  set $n+1\to n$\;
 }
 \caption{Proposed Control Policy}
\end{algorithm}

The proposed control policy is summarized in Algorithm~\ref{alg}, and Fig.~\ref{fig:con} shows a control example of the proposed policy, in which the battery follows the regulation instruction until the distance between its maximum and minimum SoC reaches $\hat{u}$.
The detailed formulation is as follows. We assume at a particular control interval $n$, $e_n$ and $r_n$ are observed, and the proposed regulation policy has the following form: $g_n(e_n,r_n)=\begin{bmatrix} c^g_n & d^g_n \end{bmatrix}$.

The control policy employs the following strategy
\begin{align}\label{eq:pol5}
    \text{If $r_n \geq 0$, } 
		c_n^g &=\min\Big\{\frac{1}{T\eta\ud{c}}(\overline{e}^g-e_n), r_n\Big\} \\
    \text{If $r_n < 0$, } 
		d_n^g &=\min\Big\{\frac{\eta\ud{d}}{T}(e_n-\underline{e}^g), r_n\Big\}
\end{align}
where $\overline{e}^g_n$ and $\underline{e}^g_n$ are the upper and lower storage energy level bound determined by the controller at the control interval $n$ for enforcing the SoC band $\hat{u}$
\begin{align}\label{eq:pol6}
    \overline{e}^g_n &= \min\{\overline{e}, e\up{min}_n + \hat{u}E\}\nonumber\\
    \underline{e}^g_n &= \max\{\underline{e}, e\up{max}_n - \hat{u}E\}\,
\end{align}
and $e\up{max}_n$, $e\up{min}_n$ is the current maximum and minimum battery storage level since the beginning of the operation, which are updated at each control step as
\begin{align}\label{eq:pol4}
    e\up{max}_n &= \max\{e\up{max}_{n-1}, e_n\}\nonumber\\
    e\up{min}_n &= \min\{e\up{min}_{n-1}, e_n\}\,.
\end{align}

\subsection{Optimality Gap to Offline Problem}

Let $(\mathbf{c}^*, \mathbf{d}^*)$ be an offline minimizer to the regulation response problem as
\begin{align}
    (\mathbf{c}^*, \mathbf{d}^*) \in \mathrm{arg}&\min_{\mathbf{c}, \mathbf{d}} J(\mathbf{c}, \mathbf{d}, \mathbf{r})\nonumber\\
    & \text{subjects to \eqref{eq:bat_soc}--\eqref{eq:bat_con4}}\label{eq:opt1}
\end{align}
and let $g(e_0, \mathbf{r})$ denote the control action of the proposed policy subjects to the initial storage energy level $e_0$ and the regulation instruction realization $\mathbf{r}$.


\begin{figure}[t]%
	\centering
	\subfloat[Regulation instruction]{
		\includegraphics[trim = 3mm 0mm 5mm 0mm, clip, width = .35\columnwidth]{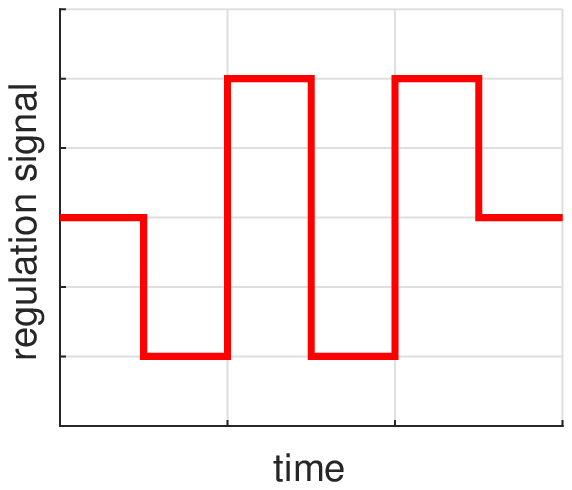}
		\label{fig:policy_A}%
	}
	\subfloat[$\theta=\pi>0$]{
		\includegraphics[trim = 3mm 0mm 5mm 0mm, clip, width = .35\columnwidth]{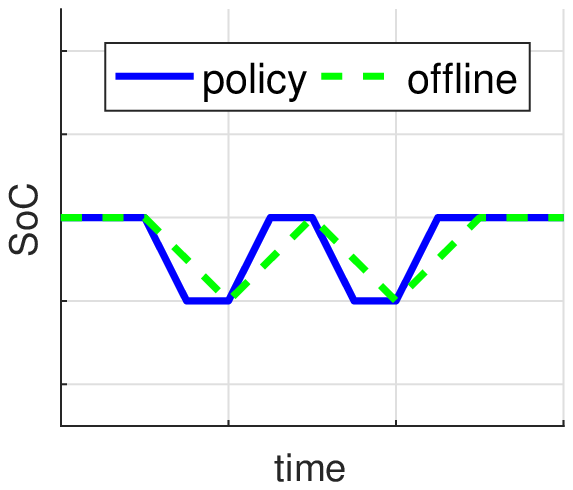}
		\label{fig:policy_B}%
	}
	\\
	\subfloat[$\theta=0, \pi>0$]{
		\includegraphics[trim = 3mm 0mm 5mm 0mm, clip, width = .35\columnwidth]{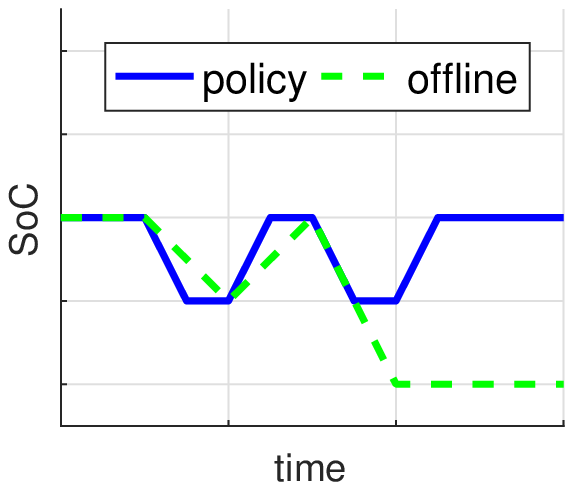}
		\label{fig:policy_C}%
	}
	\subfloat[$\theta>0, \pi=0$]{
		\includegraphics[trim = 3mm 0mm 5mm 0mm, clip, width = .35\columnwidth]{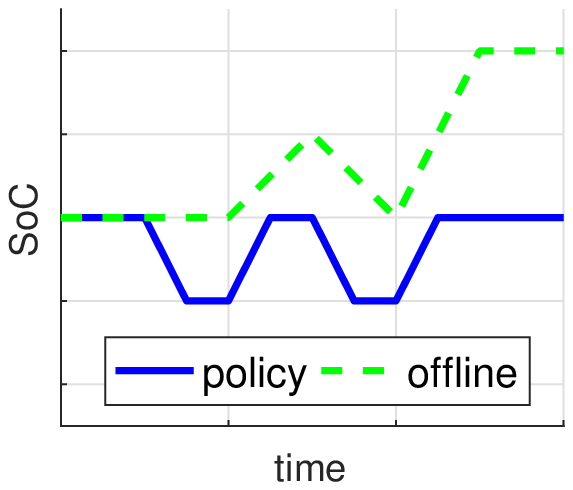}
		\label{fig:policy_D}%
	}
	\caption{Example illustration of the policy optimality under different price settings. The value of $\theta+\pi$ is the same in all cases and the round-trip efficiency is assumed to be one, so $\hat{u}$ is the same in all cases. }%
	\label{fig:policy}
\end{figure}

\begin{theorem}\label{theorem1}
Suppose the battery cycle aging stress function $\Phi(\cdot)$ is strictly convex. The proposed control strategy $g(\cdot)$ has a worst-case optimality gap (regret) $\epsilon$ that is independent of the operation time duration $TN$:
\begin{align}
	&\text{$\exists\, \epsilon \in \mathbb{R}^+$ s.t. } J(g(e_0, \mathbf{r}), \mathbf{r})-J(\mathbf{c}^*, \mathbf{d}^*, \mathbf{r}) \leq \epsilon\, \label{eq:gap}\\
    & \text{$\forall$ $e_0 \in [\underline{e}, \overline{e}]$, and  $\forall$ sequences $\{r_n\} \in [-P,P]$, $N\in\mathbb{N}$.} \nonumber
\end{align}
\end{theorem}

{The bound $\epsilon$ in Theorem~\ref{theorem1} can be explicitly characterized. To do this, we transform the objective function $J(\cdot)$ using cycle depths $(\mathbf{u},\mathbf{v},\mathbf{w})$ as control variables instead of $(\mathbf{c}, \mathbf{d})$:}
\begin{subequations}\label{eq:cyc_cost}
\begin{align}
    J\ud{u}(u) &=  ER\Phi(u) + E(\theta/\eta\ud{c}+\pi\eta\ud{d})u \label{eq:cyc_cost1}
    \\
    J\ud{v}(v) &= (1/2)ER\Phi(v) + (E/\eta\ud{c})\theta v \label{eq:cyc_cost2}
    \\
    J\ud{w}(w) &= (1/2)ER\Phi(w) + E\eta\ud{d}\pi w\, \label{eq:cyc_cost3}
\end{align}
\end{subequations}
where $J\ud{u}$ is the cost associated with a full cycle, $J\ud{v}$ for a charge half cycle, and $J\ud{w}$ for a discharge half cycle. The detailed transforming procedure is discussed in the appendix.

If we assume the cycle depth strss function $\Phi(\cdot)$ is strictly convex, then it is easy to see that \eqref{eq:pol3} is the unconstrained minimizer to \eqref{eq:cyc_cost1}. Similarly, the unconstrained minimizers of  \eqref{eq:cyc_cost2} and \eqref{eq:cyc_cost3} are:
\begin{align}
    \hat{v} = \dot{\Phi}^{-1}\Big(\frac{\theta/\eta\ud{c}}{R}\Big),\quad \hat{w} = \dot{\Phi}^{-1}\Big(\frac{\pi\eta\ud{d}}{R}\Big)\,.
\end{align}
The follow theorem offers the analytical expression for $\epsilon$
\begin{theorem}\label{theorem2}
If function $\Phi(\cdot)$ is strictly convex, then the worst-case optimality gap for the proposed policy $g(\cdot)$ as in \eqref{eq:gap} is
\begin{align}
    \epsilon = \begin{cases}
    \epsilon\ud{w} & \text{if $\pi\eta\ud{d} > \theta/\eta\ud{c}$} \\
    0 & \text{if $\pi\eta\ud{d} = \theta/\eta\ud{c}$} \\
    \epsilon\ud{v} & \text{if $\pi\eta\ud{d} < \theta/\eta\ud{c}$}
    \end{cases}\,
\end{align}
where
\begin{align}
    \epsilon\ud{w} &= J\ud{w}(\hat{u})+2J\ud{v}(\hat{u})-J\ud{w}(\hat{w})-2J\ud{v}(\hat{v}) \label{eq:gap_v}\\
    \epsilon\ud{v} &= 2J\ud{w}(\hat{u})+J\ud{v}(\hat{u})-2J\ud{w}(\hat{w})-J\ud{v}(\hat{v})\label{eq:gap_u}\,.
\end{align}
\end{theorem}

We defer the proof of this theorem to the appendix since it is somewhat technically involved.  The intuition is that battery operations consist mostly full cycles due to limited storage capacity because the battery has to be charged up before discharged, and vice versa. Enforcing $\hat{u}$--the optimal full cycle depth calculated from penalty prices and battery coefficients--ensures optimal responses in all full cycles. In cases that $\pi\eta\ud{d} = \theta/\eta\ud{c}$, $\hat{u}$ is also the optimal depth for half cycles, and the proposed policy achieves optimal control. In other cases, the optimality gap is caused by half cycles because they have different optimal depths. However, half cycles have limited occurrences in a battery operation because they are incomplete cycles~\cite{amzallag1994standardization}, so that the optimality gap is bounded as stated in Theorem~\ref{theorem2}. Fig.~\ref{fig:policy} shows some examples of the policy optimality when responding to the regulation instruction (Fig~\ref{fig:policy_A}) under different price settings. The proposed policy has the same control action in all three price settings because of the same $\hat{u}$. The policy achieves optimal control in Fig~\ref{fig:policy_B} because $\hat{u}$ is the optimal depth for all cycles. In Fig~\ref{fig:policy_C} and Fig~\ref{fig:policy_D}, half cycles have different optimal depths and the policy is only near-optimal. However, the offline result also selectively responses to instructions with a zero penalty price (charge instructions in Fig~\ref{fig:policy_C}, discharge instructions in Fig~\ref{fig:policy_D}) because it returns the battery to a shallower cycle depth so that the battery have smaller marginal operating cost in future operations.

\section{Simulation Results}\label{sec:simulation}

\begin{table*}[!htb]
    \centering
    \caption{Simulation with Random Generated Regulation Signals.}
    \begin{tabular}{l c c c c c c c c c c c}
        \hline
        \hline
         & $\theta$ & $\pi$ & $\eta$ & $N$ & $\hat{u}$ & $\epsilon$ & \multicolumn{2}{c}{ Maximum optimality gap [\$]} & \multicolumn{3}{c}{Average objective value [\$]}
        \Tstrut\\
        Case & [\$/MWh] & [\$/MWh] & [\%] & & [\%] & [\$] & Theoretical & Simple & Offline & Proposed & Simple
        \Bstrut\\
        \hline
        1 & 50 & 50 & 100 & 100 & 11.1 & \highlight{0.00} & \highlight{0.00} & \compare{183.9} & 117.4 & 117.4 & 200.2 \Tstrut\\
        2 & 100 & 100 & 100 & 100 & 21.9 & \highlight{0.00} & \highlight{0.00} & \compare{127.5} & 168.7 & 168.7 & 209.0 \\
        3 & 200 & 200 & 100 & 100 & 42.8 & \highlight{0.00} & \highlight{0.00} & \compare{$\,\,\,$47.9} & 219.4 & 219.4 & 226.7  \\
        4 & 50 & 50 & 85 & 100 & 11.2 & \highlight{0.06} & \highlight{0.06} & \compare{184.9} & 117.2 & 117.3 & 202.9 \\
        5 & 80 & 20 & 85 & 100 & 11.7 & \highlight{3.83} & \highlight{3.83} & \compare{181.4} & 108.0 & 110.7 & 198.9 \\
        6 & 20 & 80 & 85 & 100 & 10.6 & \highlight{2.19} & \highlight{2.19} & \compare{192.6} & 122.4 & 123.8 & 206.8 \\
        7 & 50 & 50 & 85 & 200 & 11.2 & \highlight{0.06} & \highlight{0.06} & \compare{408.8} & 235.6 & 235.7 & 388.3  \\
        8 & 80 & 20 & 85 & 200 & 11.7 & \highlight{3.83} & \highlight{3.83} & \compare{400.6} & 219.5 & 222.2 & 375.4  \\
        9 & 20 & 80 & 85 & 200 & 10.6 & \highlight{2.19} & \highlight{2.19} & \compare{421.4} & 247.6 & 248.9 & 401.1  \Bstrut\\
         \hline
         \hline
    \end{tabular}
    \label{tab:sim}
\end{table*}

\subsection{Simulation Setting}
We compare the proposed control policy with the offline optimal result and a simple control policy~\cite{bitar2011role}.
The maximum storage level $\overline{e}$ is set to 0.95~MWh and the minimum storage level $\underline{e}$ is set to 0.1~MWh. This assumed battery storage consists of lithium-ion battery cells that can perform 3000 cycles at 80\% cycle depth before reaching end of life, and these cells have a polynomial cycle depth stress function concluded from lab tests~\cite{laresgoiti2015modeling}:$\Phi(u) = \mathrm{(5.24\times10^{-4})}u^{2.03}$,
and the cell replacement price is set to 300~\$/kWh.

\subsection{Optimality Gap}

\begin{figure}[t]
    \centering
    \includegraphics[trim = 5mm 0mm 10mm 0mm, clip, width = .7\columnwidth]{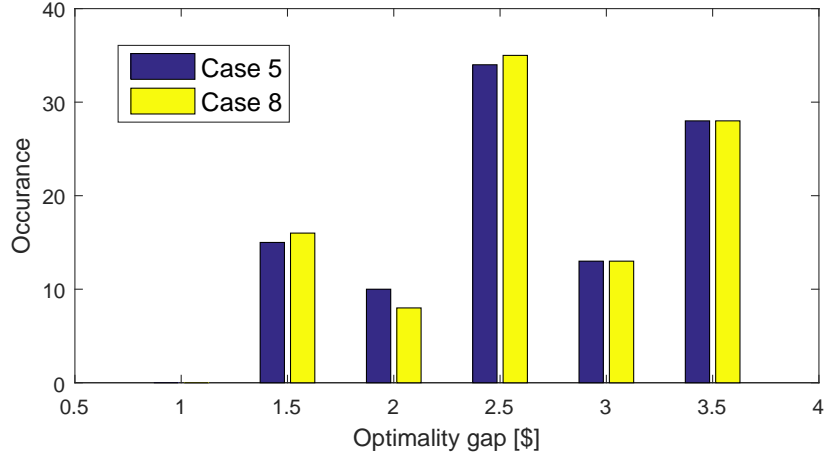}
    \caption{Distribution of optimality gaps for Case 5 and Case 8. Although Case 8 has twice the duration of Case 5, their optimality gaps are similar. This validates that the worst-case optimality gap is independent of the regulation operation duration.}
    \label{fig:gap}
\end{figure}

We simulate regulation using random generated regulation traces to exam the optimality of the proposed policy and to validate Theorem~\ref{theorem1} and \ref{theorem2}. We generate 100 regulation signal traces assuming a uniform distribution between $[-1, 1]$, and design nine test cases. Each test case has different market prices and battery round-trip efficiency $\eta=\eta\ud{d}\eta\ud{c}$. In order to demonstrate the time-invariant property of the optimality gap, Case 7 to 9 are designed to have twice the duration of Case 1 to 6 by repeating the generated regulation signal trace.

The 100 generated regulation traces are simulated using the proposed policy, the simple policy, and the offline solver for each test case. Table~\ref{tab:sim} summarizes the test results. The penalty prices, round-trip efficiency, and the number of simulation control intervals used in each test case are listed, as well as the control SoC bound $\hat{u}$ and the worst-case optimality gap $\epsilon$ that are calculated using simulation parameters. The simulation results are recorded under the maximum optimality gap and the average objective value.

\begin{figure}[t]
    \centering
    \includegraphics[trim = 5mm 0mm 10mm 0mm, clip, width = .7\columnwidth]{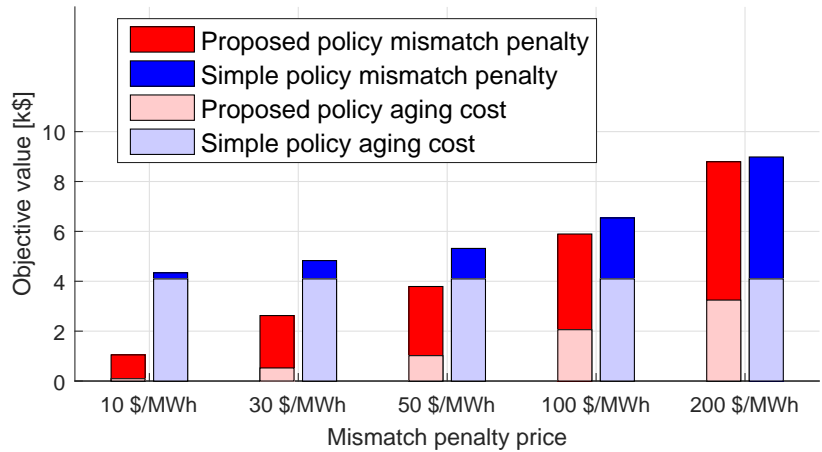}
    \caption{Regulation operating cost break-down comparison between the proposed policy and the simple policy. Although the proposed policy has higher penalties, the cost of cycle aging is significantly smaller, so it achieves better trade-offs between degradation and mismatch penalty.}
    \label{fig:pjm}
\end{figure}

This test validates Theorem~\ref{theorem2} since $\epsilon$ is exactly the same as the recorded maximum optimality gap for the proposed policy in all cases (both highlighted in pink), while the gap for the simple policy is significantly larger (highlighted in grey). In particular, the proposed policy achieves exact control results in Case 1 to 3 because $\theta/\eta\ud{c}=\pi\eta\ud{d}$, while Case 4 to 9 have non-zero gaps because the round-trip efficiency is less than ideal ($\eta<1$). We also see that as penalty prices become higher, the control band $\hat{u}$ becomes wider and the battery follows the regulation instruction more accurately. The simple policy also provides better control results at higher penalty prices.  Case 7 to 9 have the same parameter settings as to Case 4 to 6, except that the regulation signal is repeated once more time. The proposed policy achieves the same worst-case optimality gap in the two operation duration settings, while the average objective values are approximately doubled as shown in Fig.~\ref{fig:gap}.

\subsection{Simulation using Realistic Regulation Signal}

In this simulation we compare the proposed policy with the simple policy using the regulation signal trace published by PJM Interconnection~\cite{pjm_signal}. The control time interval for this signal is 2 seconds and the duration is 4 weeks. We do not use the offline result for comparison in this case because this problem is beyond the solvability of the implemented numerical solver.

We repeat the simulation using different penalty prices. We let $\theta=\pi$ in each test case and set the round-trip efficiency to $85\%$. Fig.~\ref{fig:pjm} summarizes the simulation results in the form of regulation operating cost versus penalty prices, the cycle aging cost and the regulation mismatch penalty are listed for each policy. Because the simple control does not consider market prices, its control actions are the same in all price scenarios, and the penalty increases linearly with the penalty price. The proposed policy causes significantly smaller cycle aging cost, and have better control results. As the penalty price increases, the gap between the two policies becomes smaller since $\hat{u}$ becomes closer to 100\%. 

\section{Conclusion}\label{sec:con}
In this paper, we proposed an online policy for a battery owner to provide frequency regulation in a pay-for-performance market. It considers the cycle aging mechanism of electrochemical battery cells, and is adaptive to changing market prices. We have shown that the proposed policy has a bounded regret that is dependent of operation durations, and achieves exact control result under certain market scenarios. The proposed policy applies to all battery applications that has constant prices over a specific period and the battery is dispatched in stochastic manners, such as using co-located battery storage to smooth wind farm power productions, or using behind-the-meter batteries to improve the penetration of roof-top photovoltaic generations. In our future work, we will investigate how to incorporate the policy into problems such as optimal battery contracting, and optimal battery sizing.

\bibliographystyle{IEEEtran}	
\bibliography{IEEEabrv,literature}		

\appendix

\section{Proof for Theorem~\ref{theorem1}}


\subsection{Model Reformulation}
We rewrite problem \eqref{eq:opt1} when $\pi\geq0$, $\theta\geq 0$ as 
\begin{subequations}\label{eq:opt2}
\begin{align}
    (\mathbf{c}^*, \mathbf{d}^*) \in \mathrm{arg} &\min_{\mathbf{c}, \mathbf{d}} J\ud{cyc}(\mathbf{c}, \mathbf{d}) - T\sum_{n=1}^N\big[\theta c_n + \pi d_n\big]\label{eq:the1_0}\\
    &\text{subject to \eqref{eq:bat_soc}, \eqref{eq:bat_con1}, and} \nonumber\\
    & 0\leq c_n \leq [r_n]^+  \label{eq:the1_1}\\
    & 0\leq d_n \leq [-r_n]^+ \label{eq:the1_2}
\end{align}
\end{subequations}
by observing that a battery's actions would never exceed the regulation signals.

We utilize the rainflow algorithm to transform the problem into a cycle-based form. The rainflow method maps the entire operation uniquely to cycles, the sum of all charge and discharge power can be represented as the sum of cycle depths as (recall that a full cycle has symmetric depth for charge and discharge)
\begin{align}
    \textstyle\sum_{i=1}^{|\mathbf{u}|}u_i + \textstyle\sum_{i=1}^{|\mathbf{v}|}v_i &= \frac{T\eta\ud{c}}{E}\textstyle\sum_{n=1}^N c_n
    \label{eq:rf1}\\
    \textstyle\sum_{i=1}^{|\mathbf{u}|}u_j + \textstyle\sum_{i=1}^{|\mathbf{w}|}w_i &= \frac{T}{\eta\ud{d}E}\textstyle\sum_{n=1}^N d_n\,.\label{eq:rf2}
\end{align}
We substitute \eqref{eq:rf1} and \eqref{eq:rf2} into the reformulated objective function \eqref{eq:the1_0} to replace $c_n$ and $d_n$ with cycle depths
\begin{align}
    &J\ud{cyc}(\mathbf{c}, \mathbf{d}) + J\ud{reg}(\mathbf{c}, \mathbf{d}, \mathbf{r}) =\nonumber\\ &\textstyle\sum_{i=1}^{|\mathbf{u}|}J\ud{u}(u_i) + \textstyle\sum_{i=1}^{|\mathbf{v}|}J\ud{v}(v_i) + \textstyle\sum_{i=1}^{|\mathbf{w}|}J\ud{w}(w_i)\,.
\end{align}

\subsection{Proof for Theorem~\ref{theorem2}}

The following lemmas support the proof for Theorem~\ref{theorem2}
\begin{lemma}\label{lemma1}
Assume a minimizer $(\mathbf{c}^*, \mathbf{d}^*)$ to problem \eqref{eq:opt1} has the cycle counting results $(\mathbf{u}^*, \mathbf{v}^*, \mathbf{w}^*)$. Then the depth of each cycle in this result either reaches the optimal cycle depth or bounded by the operation constraints \eqref{eq:bat_con1}, \eqref{eq:the1_1}, \eqref{eq:the1_2} as
\begin{subequations}
\begin{align}
    u_i^* &= \min(\hat{u}, \overline{u}_i) \\
    v_i^* &= \min(\hat{v}, \overline{v}_i) \\
    w_i^* &= \min(\hat{w}, \overline{w}_i)
\end{align}
\end{subequations}
where $\overline{u}_i$, $\overline{v}_i$, $\overline{w}_i$ denote constraint bounds.
\end{lemma}

\begin{lemma}\label{lemma2}
A cycle depth in the control action of $g(\cdot)$ either reaches the depth of $\hat{u}$ or is bounded by the operation constraints.
\end{lemma}

\begin{lemma}\label{lemma3}
There exists one and only one half cycle with the largest depth in a rainflow residue profile. Other half cycles are in strictly decreasing order either to the left- or to the right-hand side direction of this largest half cycle.
\begin{figure}[ht]
    \centering
    \includegraphics[trim = 5mm 0mm 10mm 0mm, clip, width = .7\columnwidth]{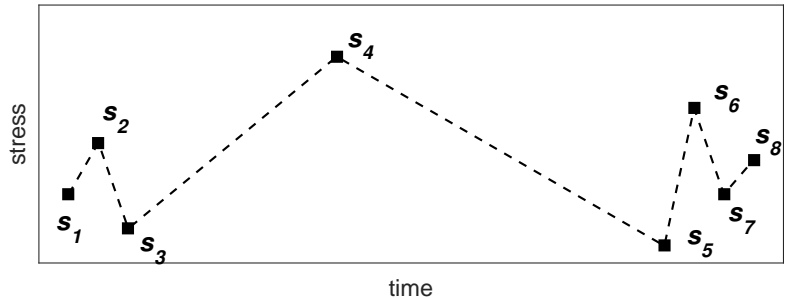}
    \caption{Illustration for Lemma~\ref{lemma3}. The largest half cycle is between $s_4$ and $s_5$, other half cycles are in strictly decreasing order either to the left- or to the right-hand side direction of this largest half cycle.}
    \label{fig:lemma3}
\end{figure}
\end{lemma}

It is easy to see now from Lemma~\ref{lemma1} and Lemma~\ref{lemma2} that the proposed control policy achieves optimal control result for all full cycles, and the optimality gap is caused by half cycle results. Consider the following relationship in a rainflow residue profile as in Lemma~\ref{lemma3} assuming the largest half cycle is in the discharging direction
\begin{align}\label{the3:eq2}
    \dotsc < w^*_{j-1} < v^*_{j-1} <w^*_j>v^*_j >w^*_{j+1}>\dotsc
\end{align}
and substitute Lemma~\ref{lemma2} into \eqref{the3:eq2}
\begin{align}
   \dotsc \min\{\hat{v}, \overline{v}_j\} <\min\{\hat{w}, \overline{w}_j\}>\min\{\hat{v}, \overline{v}_{j+1}\}\dotsc
\end{align}
It is easy to see now that if $\hat{w}>\hat{v}$, then the largest possible value for $ w^*_j$ is $\hat{w}$,
and the largest possible value for $v^*_j$ and $v^*_{j-1}$ is $\hat{v}$, the rest half cycles in \eqref{the3:eq2} must have depths smaller than $\hat{v}$, which indicates that their depths are bounded by operation. If $\hat{v}>\hat{w}$, then the largest possible value for $ w^*_j$ is $\hat{w}$, and the rest half cycles must have depths smaller than $\hat{w}$. We repeat this analysis for cases that $v^*_j$ is the largest cycle, and summarize the half cycle conditions in Table~\ref{tab:half_cycle}
\begin{table}[ht]
    \centering
    \caption{Summarizing Half Cycle Depth Conditions}
    \begin{tabular}{l c c}
        \hline
        \hline
        &  $\hat{w}>\hat{v}$ & $\hat{w}<\hat{v}$ \Tstrut\Bstrut\\
        \hline
        Half cycles of depth $\hat{w}$ & At most one & At most two\Tstrut\Bstrut\\
        Half cycles of depth $\hat{v}$ & At most two & At most one\Tstrut\Bstrut\\
        Rest half cycles & must be $<\hat{v}$ & must be $<\hat{w}$\Tstrut\Bstrut\\
        \hline
        \hline
    \end{tabular}
    \label{tab:half_cycle}
\end{table}
Hence, the worst-case optimality gap is caused by that some half cycles have depth $\hat{u}$ or $\hat{w}$, while the control policy enforces $\hat{u}$ as the depth of all cycles unbounded by operation. The gap in Theorem~\ref{theorem2} is therefore calculated using half cycle depth conditions in Table~\ref{tab:half_cycle}.

\noindent \emph{Proof for Lemma~\ref{lemma1}:}
This property is trivial because cycles are linear combinations of charge and discharge power, and constraints \eqref{eq:bat_con1}, \eqref{eq:the1_1}, \eqref{eq:the1_2} can be transformed into linear constrains with respect to cycle depths. Hence the transformed cycle-based problem also has a convex objective function with linear constraints. Although exact formulations of the transformed constraints are complicated to express, we use $\overline{u}_i$, $\overline{v}_i$, and $\overline{w}_i$ to denote these binds, which are sufficient for the proof for Theorem~\ref{theorem2}.\\

\noindent \emph{Proof for Lemma~\ref{lemma2}:}
\begin{figure}[ht]
    \centering
    \includegraphics[trim = 5mm 0mm 10mm 0mm, clip, width = .7\columnwidth]{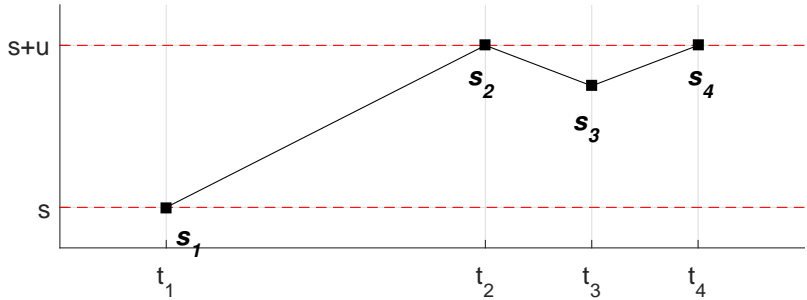}
    \caption{Illustration for Lemma~\ref{lemma2}.}
    \label{fig:lemma2}
\end{figure}
The rainflow method always identify the largest cycle as between the minimum and the maximum SoC point. In the proposed policy, any operation that goes outside the defined operation zone will cause the largest cycle depth to change instead of the depth of the cycle it was previous in. For example, in Fig.~\ref{fig:lemma2} the maximum cycle is between SoC $s$ and $s+u$, and the battery is at time $t_4$. If the battery continue to charge and the SoC goes about $s+u$, then this operation will increase the largest cycle depth instead of the shallower cycles assoicated with extrema $s_2$, $s_3$ and $s_4$.\\

\noindent \emph{Proof for Lemma~\ref{lemma3}:}
Because the rainflow method identifies a cycle from extrema distances if $\Delta s_{i-1}\geq \Delta s_{i} \leq \Delta s_{i+1}$, then all extrema in the rainflow residue must satisfy either $\Delta s_{i-1}< \Delta s_{i}  < s_{i+1}$ or $\Delta s_{i-1}< \Delta s_{i}  > s_{i+1}$ or $\Delta s_{i-1}> \Delta s_{i} > s_{i+1}$, which proofs this lemma.

\end{document}